\documentclass[10pt]{article}
\begin{document}
\title{ {\bf  On Potentially $(K_5-C_4)$-graphic
Sequences}
\thanks{  Project Supported by NNSF of China(10271105), NSF of Fujian(Z0511034),
Science and Technology Project of Fujian, Project of Fujian
Education Department and Project of Zhangzhou Teachers College.}}
\author{{ Lili Hu , Chunhui Lai}\\
{\small Department of Mathematics, Zhangzhou Teachers College,}
\\{\small Zhangzhou, Fujian 363000,
 P. R. of CHINA.}\\{\small  jackey2591924@163.com ( Lili Hu)}
 \\{\small   zjlaichu@public.zzptt.fj.cn(Chunhui
 Lai, Corresponding author)}
}
\date{}
\maketitle
\begin{center}
\begin{minipage}{4.1in}
\vskip 0.1in
\begin{center}{\bf Abstract}\end{center}
 { In this paper, we characterize
 the potentially $(K_5-C_4)$-graphic sequences where
 $K_5-C_4$ is the graph obtained from $K_5$ by removing four
 edges of a 4 cycle $C_4$. This characterization implies a
 theorem due to  Lai [6].
} \par
\par
 {\bf Key words:} graph; degree sequence; potentially $(K_5-C_4)$-graphic
sequences\par
  {\bf AMS Subject Classifications:} 05C07\par
\end{minipage}
\end{center}
 \par
 \section{Introduction}
\par

   An $n$-term nonincreasing nonnegative integer sequence
 $\pi=(d_1,d_2,\cdots,d_n)$ is said to be graphic if it is the
 degree sequence of a simple graph $G$ of order $n$; such a
 graph $G$ is referred to as a realization of $\pi$. We denote by
 $\sigma(\pi)$ the sum of all the terms of $\pi$. $K_n$ is the
 complete graph on $n$ vertices. $C_n$ is the cycle of length $n$. $K_n-C_4$ is the graph obtained from $K_n$ by removing
 4
 edges of a 4 cycle $C_4$. Let $H$
 be a simple graph. A graphic sequence $\pi$ is said to be
 potentially $H$-graphic if it has a realization $G$
 containing $H$ as a subgraph.

\par

  Given a graph $H$, what is the maximum number of edges of a graph
with $n$ vertices not containing $H$ as a subgraph? This number is
denoted $ex(n,H)$, and is known as the Tur\'{a}n number. This
problem was proposed for $H = C_4$ by Erd\"os [1] in 1938 and
generalized by Tur\'{a}n [16]. In terms of graphic sequences, the
number $2ex(n,H)+2$ is the minimum even integer $l$ such that every
$n$-term graphical sequence $\pi$ with $\sigma (\pi)\geq l $ is
forcibly $H$-graphical. In [3], Gould, Jacobson and Lehel considered
the following variation of the classical Tur\'{a}n-type extremal
problems: determine the smallest even integer $\sigma(H,n)$ such
that every n-term positive graphic sequence
$\pi=(d_1,d_2,\cdots,d_n)$ with $\sigma(\pi)\geq \sigma(H,n)$ has a
realization $G$ containing $H$ as a subgraph. They proved that
$\sigma(pK_2, n)=(p-1)(2n-p)+2$ for $p\ge 2$; $\sigma(C_4,
n)=2[{{3n-1}\over 2}]$ for $n\ge 4$.  In [5,6], Lai determined the
values $\sigma(K_4-e,n)$ for $n\geq4$ and $\sigma(K_5-C_4,n)$ for
$n\geq5$. Yin, Li, and Mao [14] determined the values
$\sigma(K_{r+1}-e,n)$ for $r\geq3$ and $r+1\leq n\leq 2r$ and
$\sigma(K_5-e,n)$ for $n\geq5$. Recently, Yin and Li [15] determined
$\sigma(K_{r+1}-e,n)$.
 Erd\"os,\ Jacobson and Lehel [2] showed that $\sigma(K_k, n)\ge
(k-2)(2n-k+1)+2$ and conjectured that the equality holds. They
proved the conjecture is true for $k=3$ and $n\geq6$, i.e.,
$\sigma(K_3,n)=2n$ for $n\geq6$. The conjecture was confirmed in
[3], [7], [8], [9] and [10].

 \par

   Motivated by the above problems, we consider the following
 problem: given a graph $H$, characterize the potentially
 $H$-graphic sequences without zero terms. In [11], Luo characterized the potentially
 $C_k$-graphic sequences for each $k=3,4,5$. Recently, Luo and Warner [12] characterized the potentially
 $K_4$-graphic sequences. In [13], Eschen and Niu characterized the potentially
 $(K_4-e)$-graphic sequences.
\par

   In this paper, we characterize the potentially $(K_5-C_4)$-graphic
 sequences without zero terms. This characterization implies a
 theorem due to  Lai [6].

\par

\section{Preparations}\par
   Let $\pi=(d_1,d_2,\cdots,d_n)$ be a nonincreasing positive
 integer sequence. Then
 $\pi^\prime=(d_1-1,d_2-1,\cdots,d_{d_{n}}-1,d_{d_{n}+1},\cdots,d_{n-1})$
 is the residual sequence obtained by laying off $d_n$ from $\pi$.
 We denote the nonincreasing sequence $\pi^\prime$ by
 $(d_1^\prime,d_2^\prime,\cdots,d_{n-1}^\prime)$. From here on,
 denote $\pi^\prime$ the residual sequence obtained by laying off
 $d_n$ from $\pi$ and all the graphic sequences have no zero terms. In order to prove
 our main result, we need the following results.\par
 \par

    \par
    {\bf Theorem 2.1} [3] If $\pi=(d_1,d_2,\cdots,d_n)$ is a graphic
 sequence with a realization $G$ containing $H$ as a subgraph,
 then there exists a realization $G^\prime$ of $\pi$ containing $H$ as a
 subgraph so that the vertices of $H$ have the largest degrees of
 $\pi$.\par
 \par

    \par
    The following corollary is obvious.\par

    \par
    \par
    {\bf Corollary 2.2}   Let $H$ be a simple graph. If $\pi^\prime$ is
 potentially $H$-graphic, then $\pi$ is
 potentially $H$-graphic.
    \par
    We will use Corollary 2.2 repeatedly in the proofs of our main
results.\par
\par

\par
   {\bf Lemma 2.3} (Kleitman and Wang [4])  $\pi$ is
graphic if and only if $\pi^\prime$ is graphic.
\par
\par

    \par
\section{  Potentially $(K_5-C_4)$-graphic sequences} \par

  Our main result is as follows:\par

\par
  \textbf{Theorem 3.1}  Let $\pi=(d_1,d_2,\cdots,d_n)$ be a graphic sequence with
 $n\geq5$. Then $\pi$ is potentially $(K_5-C_4)$-graphic if and
 only if the following conditions hold:

  \par
 $(1)$ $d_1\geq4$.
  \par
  $(2)$ $d_5\geq2$.

  \par
   $(3)$ $\pi\neq((n-2)^2,2^{n-2})$ for $n\geq6$, where the symbol
   $x^y$ stands for $y$ consecutive terms $x$.

  \par
   $(4)$ $\pi\neq(n-k,k+i,2^i,1^{n-i-2})$ where
   $i=3,4,\cdots,n-2k$ and $k=1,2,\cdots,[\frac{n-1}{2}]-1$.

  \par
   $(5)$ If $n=6$, then $\pi\neq (4,2^5)$.

  \par
   $(6)$ If $n=7$, then $\pi\neq(4,2^6)$.
    \par
    Proof: First we assume that $\pi$ is potentially
$(K_5-C_4)$-graphic. In this case the necessary conditions $(1)$
and $(2)$ are obvious.  we are going to prove the conditions
$(3)-(6)$ by way of contradiction.

  \par
    If $\pi=((n-2)^2,2^{n-2})$ where $n\geq6$ is potentially
$(K_5-C_4)$-graphic, then according to theorem 2.1, there exists a
realization $G$ of $\pi$ containing $K_5-C_4$ as a subgraph so
that the vertices of $K_5-C_4$ have the largest degrees of $\pi$.
Then the sequence $\pi^*=(n-4,n-6,2^{n-5})$ obtained from
$G-(K_5-C_4)$ must be graphic and there must be no edge between
two vertices with degree $n-4$ and $n-6$ for the realization of
$\pi^*$, which is impossible. Thus, $\pi=((n-2)^2,2^{n-2})$ where
$n\geq6$ is not potentially $(K_5-C_4)$-graphic. Hence, $(3)$
holds.
  \par
  If $\pi=(n-k,k+i,2^i,1^{n-i-2})$ where $i=3,4,\cdots,n-2k$ and $k=1,2,\cdots,[\frac{n-1}{2}]-1$ is potentially
$(K_5-C_4)$-graphic, then according to theorem 2.1, there exists a
realization $G$ of $\pi$ containing $K_5-C_4$ as a subgraph so
that the vertices of $K_5-C_4$ have the largest degrees of $\pi$.
Then the sequence $\pi^*=(n-k-4,k+i-2,2^{i-3},1^{n-i-2})$ obtained
from $G-(K_5-C_4)$ must be graphic and there must be no edge
between two vertices with degree $n-k-4$ and $k+i-2$ for the
realization of $\pi^*$. Thus, $\pi^*$ satisfies:
$(n-k-4)+(k+i-2)\leq2(i-3)+(n-i-2)$, that is, $0\leq(-2)$, which
is a contradiction. Hence, $(4)$ holds.
\par
  If  $\pi=(4,2^5)$ is potentially $(K_5-C_4)$-graphic, then
according to theorem 2.1, there exists a realization $G$ of $\pi$
containing $K_5-C_4$ as a subgraph so that the vertices of
$K_5-C_4$ have the largest degrees of $\pi$. Then the sequence
$\pi^*=(2)$ obtained from $G-(K_5-C_4)$ must be the degree
sequence of a simple graph, which is a contradiction. Thus,
$\pi=(4,2^5)$ is not potentially $(K_5-C_4)$-graphic. Hence, $(5)$
holds.

  \par
  If  $\pi=(4,2^6)$ is potentially
$(K_5-C_4)$-graphic, then according to theorem 2.1, there exists a
realization $G$ of $\pi$ containing $K_5-C_4$ as a subgraph so
that the vertices of $K_5-C_4$ have the largest degrees of $\pi$.
Then the sequence $\pi^*=(2^2)$ obtained from $G-(K_5-C_4)$ must
be the degree sequence of a simple graph, which is a
contradiction. Thus, $\pi=(4,2^6)$ is not potentially
$(K_5-C_4)$-graphic. Hence, $(6)$ holds.
  \par
  Now we prove the sufficient condition. Suppose the graphic
sequence $\pi$ satisfies the conditions $(1)-(6)$. Our proof is by
induction on $n$.
\par
  First we prove the sufficient condition for $n=5$. Since
$\pi\neq(4^2,2^3)$, then $\pi$ is one of the following sequences:
\par
   $(4^5)$, $(4^3,3^2)$, $(4^2,3^2,2)$, $(4,3^4)$,
   $(4,3^2,2^2)$, $(4,2^4)$. It is easy to see that they are all potentially $(K_5-C_4)$-graphic. Therefore,
$\pi$ is potentially $(K_5-C_4)$-graphic for $n=5$.
  \par
   We now suppose that the sufficient condition holds for
$(n-1)\geq5$. We will prove that it holds for $n$. Let
$\pi=(d_1,d_2,\cdots,d_n)$ be a graphic sequence with $n$ terms
that satisfies the conditions $(1)-(6)$. We only need to show that
$\pi$ is potentially $(K_5-C_4)$-graphic. If $\pi^\prime$
satisfies the assumption, then $\pi^\prime$ is potentially
$(K_5-C_4)$-graphic by the induction hypothesis. Therefore, $\pi$
is potentially $(K_5-C_4)$-graphic by Corollary 2.2. Thus, we
consider the following cases:

  \par
  \textbf{Case 1:} If $\pi^\prime=(4,2^5)$, then $\pi=(5,3,2^5)$ or $\pi=(5,2^5,1)$. It
is easy to see that both of them are potentially
$(K_5-C_4)$-graphic.

  \par
  \textbf{Case 2:} If $\pi^\prime=(4,2^6)$, then
$\pi=(5,3,2^6)$ or $\pi=(5,2^6,1)$. It is easy to see that both of
them are potentially $(K_5-C_4)$-graphic.

    \par
    \textbf{Case 3:} $\pi^\prime=((n-3)^2,2^{n-3})$ where
 $n-1\geq6$.

\par
  If $d_n=2$, then $\pi=((n-2)^2,2^{n-2})$, which is contradict to condition(3).

\par

    If  $d_n=1$, then $\pi=(n-2,n-3,2^{n-3},1)$. We are going to
prove that $\pi$ is potentially $(K_5-C_4)$-graphic. First we show
it is true for $n=6$. In this case, $\pi=(4,3,2^3,1)$. It is easy
to see that $\pi$ is potentially $(K_5-C_4)$-graphic. Now we prove
that $\pi$ is potentially $(K_5-C_4)$-graphic for $n\geq7$. It is
enough to show $\pi_1=(n-5,n-6,2^{n-6},1)$ is graphic and there
exist no edge between two vertices with degree $n-5$ and $n-6$ for
the realization of $\pi_1$. Hence it is enough to show
$\pi_2=(n-6,1^{n-6})$ is graphic. Clearly, $\pi_2$ has a
realization consisting of $n-6$ edges and these edges have only
one vertex in common.
  \par
   Thus, $\pi=(n-2, n-3,2^{n-3},1)$ is potentially $(K_5-C_4)$-graphic for $n\geq6$.

    \par
    \textbf{Case 4:} $\pi^\prime=(n-1-k,k+i,2^i,1^{n-i-3})$ where
  $i=3,4,\cdots,n-1-2k$ and $k=1,2,\cdots,[\frac{n-2}{2}]-1$.

    \par
  If $d_n=2$, then $n-i-3=0$ and $\pi=(n-k,k+i+1,2^{i+1})$, which is contradict to condition(4).

  \par
  If  $d_n=1$, then $\pi=(n-k^\prime,k^\prime+i,2^i,1^{n-i-2})$, which is contradict to condition(4) .
\par

 \textbf{Case 5:} $d_n\geq4$. In this case, $\pi^\prime$ satisfies the conditions
$(1)-(6)$. Thus, $\pi^\prime$ is potentially $(K_5-C_4)$-graphic.
Therefore, $\pi$ is potentially $(K_5-C_4)$-graphic by Corollary
2.2.

\par

 \textbf{Case 6:} $d_n=3$.

  \par
   If $d_1\geq5$, then $\pi^\prime$ satisfies the conditions
$(1)-(6)$. Thus, $\pi^\prime$ is potentially $(K_5-C_4)$-graphic.
Therefore, $\pi$ is potentially $(K_5-C_4)$-graphic by Corollary
2.2.

   \par
   If $d_1=4$, there are two subcases: $d_4=4$ and $d_4=3$.

 \par
   Subcase 1: $d_4=4$. In this case, $d_1=d_2=d_3=d_4=4$.
Obviously, $\pi^\prime$ satisfies the conditions $(1)-(6)$. Thus,
$\pi^\prime$ is potentially $(K_5-C_4)$-graphic. Therefore, $\pi$
is potentially $(K_5-C_4)$-graphic by Corollary 2.2.

 \par
   Subcase 2: $d_4=3$.

    \par
   Subcase 2.1: $d_3=4$. Then $\pi=(4^3,3^{n-3})$. Since $\sigma(\pi)$ is even, $n$ must be odd. We are going to
prove that $\pi$ is potentially $(K_5-C_4)$-graphic. It is easy to
see that $\pi=(4^3,3^4)$ is potentially $(K_5-C_4)$-graphic. If
$n\geq9$, then $(4^3,3^{n-3})$ has a realization containing a
$K_5-C_4$ (see Figure 1).
    \par
    Thus, $\pi=(4^3,3^{n-3})$ where $n$ is odd is potentially
$(K_5-C_4)$-graphic.

\par
   Subcase 2.2: $d_3=3$.

    \par
   If $d_2=4$, then $\pi=(4^2,3^{n-2})$. Since $\sigma(\pi)$ is even,  $n$ must be even. We are going to
prove that $\pi$ is potentially $(K_5-C_4)$-graphic. It is easy to
see that $\pi=(4^2,3^4)$ and $\pi=(4^2,3^6)$ are potentially
$(K_5-C_4)$-graphic. If $n\geq10$, then $(4^2,3^{n-2})$ has a
realization containing a $K_5-C_4$ (see Figure 2).

 \par
  Thus, $\pi=(4^2,3^{n-2})$ where $n$ is even is potentially
$(K_5-C_4)$-graphic.

    \par
    If $d_2=3$, then $\pi=(4,3^{n-1})$. Since $\sigma(\pi)$ is even, $n$ must be odd . We are going to
prove that $\pi$ is potentially $(K_5-C_4)$-graphic. It is easy to
see that $\pi=(4,3^6)$ is potentially $(K_5-C_4)$-graphic. If
$n\geq9$, then $(4,3^{n-1})$ has a realization containing a
$K_5-C_4$ (see Figure 3).

    \par
      Thus, $\pi=(4,3^{n-1})$ where $n$ is odd is potentially
$(K_5-C_4)$-graphic.

\par
 \textbf{Case 7:} $d_n=2$ and $\pi^\prime\neq((n-3)^2,2^{n-3})$ where $n-1\geq6$,
$\pi^\prime\neq(n-1-k,k+i,2^i,1^{n-i-3})$ where
 $i=3,4,\cdots,n-1-2k$ and $k=1,2,\cdots,[\frac{n-2}{2}]-1$. $\pi^\prime\neq(4,2^5)$, $\pi^\prime\neq(4,2^6)$.
\par

 If  $d_1\geq5$, then $\pi^\prime$ satisfies the conditions
$(1)-(6)$. Thus, $\pi^\prime$ is potentially $(K_5-C_4)$-graphic.
Therefore, $\pi$ is potentially $(K_5-C_4)$-graphic by Corollary
2.2.
\par
  If $d_1=4$, there are three subcases: $d_2=4$, $d_2=3$ and
$d_2=2$.

\par

 Subcase 1: $d_2=4$.
 \par
  If $d_3=4$, then $\pi^\prime$ satisfies the conditions
$(1)-(6)$. Thus, $\pi^\prime$ is potentially $(K_5-C_4)$-graphic.
Therefore, $\pi$ is potentially $(K_5-C_4)$-graphic by Corollary
2.2.
\par
  If $d_3=3$, then $\pi=(4^2, 3^a, 2^{n-2-a})$ where $a\geq1$ and $n-2-a\geq1$. Since $\sigma(\pi)$ is even, $a$ must be even.
We are going to prove that $\pi$ is potentially
$(K_5-C_4)$-graphic.
\par
  First, we consider $\pi=(4^2,3^2,2^{n-4})$. It is easy to see
that $\pi=(4^2,3^2,2^{2})$ and $\pi=(4^2,3^2,2^{3})$ are
potentially $(K_5-C_4)$-graphic. If $n\geq 8$, then
$(4^2,3^2,2^{n-4})$ has a realization containing a $K_5-C_4$ (see
Figure 4). Thus, we are done.
\par
  Then we consider $\pi=(4^2,3^a,2^{n-2-a})$ where $a\geq4$ and $n-2-a
  \geq1$. It is easy to
see that $\pi=(4^2,3^4,2)$ and $\pi=(4^2,3^4,2^2)$ are potentially
$(K_5-C_4)$-graphic. If $a=4$ and $n\geq9$, then
$(4^2,3^4,2^{n-6})$ has a realization containing a $K_5-C_4$ (see
Figure 5). If $a\geq6$, then $(4^2,3^a,2^{n-2-a})$ has a
realization containing a $K_5-C_4$ (see Figure 6).
\par
  If $d_3=2$, then $\pi=(4^2, 2^{n-2})$. Since $\pi\neq(4^2,2^4)$, we must have $n\geq7$. We are going to prove that $\pi$ is potentially
$(K_5-C_4)$-graphic. It is enough to show $\pi_1=(2^{n-4})$ is
graphic. Clearly, $C_{n-4}$ is a realization of $\pi_1$. Thus, we
are done.

\par
  Subcase 2:  $d_2=3$. Then $\pi=(4,3^a,2^{n-1-a})$ where $a\geq1$ and $n-1-a\geq1$. Since $\sigma(\pi)$ is even, $a$ must be even.
We are going to prove that $\pi$ is potentially
$(K_5-C_4)$-graphic.

\par
  First, we consider $\pi=(4,3^2,2^{n-3})$. It is enough to show
$\pi_1=(2^{n-5},1^2)$ is graphic. Clearly, $\pi_1$ is graphic.
Thus, $\pi=(4,3^2,2^{n-3})$ is potentially $(K_5-C_4)$-graphic.
\par
  Second, we consider $\pi=(4,3^4,2^{n-5})$. It is easy to see
that $\pi=(4,3^4,2)$ and $\pi=(4,3^4,2^{2})$ are potentially
$(K_5-C_4)$-graphic. If $n\geq 8$, then $(4,3^4,2^{n-5})$ has a
realization containing a $K_5-C_4$ (see Figure 7). Thus, we are
done.
\par
  Then we consider $\pi=(4,3^a,2^{n-1-a})$ where $a\geq6$ and $n-1-a \geq 1$.
 It is easy to see that $\pi=(4,3^6,2)$ is potentially
 $(K_5-C_4)$-graphic. If $a=6$ and $n\geq9$, then
$(4,3^6,2^{n-7})$ has a realization containing a $K_5-C_4$ (see
Figure 8). If $a\geq8$ and $n-1-a=1$, then $(4,3^a,2)$ has a
realization containing a $K_5-C_4$ (see Figure 9). If $a\geq8$ and
$n-1-a\geq2$, then $(4,3^a,2^{n-1-a})$ has a realization
containing a $K_5-C_4$ (see Figure 10).
 Thus, we
are done.

\par
  Subcase 3: $d_2=2$. Then $\pi=(4,2^{n-1})$. Since $\pi\neq (4,2^5)$ and $\pi\neq (4,2^6)$, we must have $n\geq8$. We are going to prove that $\pi$ is potentially $(K_5-C_4)$-graphic.
It is enough to show $\pi_1=(2^{n-5})$ where $n\geq8$ is graphic.
Obviously, $C_{n-5}$ is a realization of $\pi_1$. Thus,
$\pi=(4,2^{n-1})$ is potentially $(K_5-C_4)$-graphic.

\par
  \textbf{Case 8:} $d_n=1$ and $\pi^\prime\neq((n-3)^2,2^{n-3})$, $\pi^\prime\neq(n-1-k,k+i,2^i,1^{n-i-3})$ where
 $i=3,4,\cdots,n-1-2k$ and $k=1,2,\cdots,[\frac{n-2}{2}]-1$.
 $\pi^\prime\neq(4,2^5)$, $\pi^\prime\neq(4,2^6)$.
\par
  If  $d_1\geq5$, then $\pi^\prime$ satisfies the conditions
$(1)-(6)$. Thus, $\pi^\prime$ is potentially $(K_5-C_4)$-graphic.
Therefore, $\pi$ is potentially $(K_5-C_4)$-graphic by Corollary
2.2.
\par
  If $d_1=4$, there are three subcases: $d_2=4$, $d_2=3$ and
$d_2=2$.
\par
  Subcase 1: $d_2=4$. In this case, $\pi^\prime$ satisfies the conditions
$(1)-(6)$. Thus, $\pi^\prime$ is potentially $(K_5-C_4)$-graphic.
Therefore, $\pi$ is potentially $(K_5-C_4)$-graphic by Corollary
2.2.
\par
  Subcase 2: $d_2=3$. Then $\pi=(4,3^a,2^b,1^{n-1-a-b})$ where $a\geq1$,  $a+b\geq4$ and $n-1-a-b \geq 1$. Since $\sigma(\pi)$ is even,
$n-1-b$ must be even. We are going to prove that $\pi$ is
potentially $(K_5-C_4)$-graphic.
\par
  Subcase 2.1: $a=1$. Then $\pi=(4,3,2^b,1^{n-2-b})$. It is enough to show
$\pi_1=(2^{b-3},1^{n-1-b})$ is graphic.  Clearly, $\pi_1$ is
graphic. Thus, we are done.

\par
  Subcase 2.2: $a=2$. Then $\pi=(4,3^2,2^b,1^{n-3-b})$. It is enough to show
$\pi_1=(2^{b-2},1^{n-1-b})$ is graphic.  Clearly, $\pi_1$ is
graphic. Thus, we are done.
\par
  Subcase 2.3: $a=3$. Then $\pi=(4,3^3,2^b,1^{n-4-b})$. First, we consider $\pi=(4,3^3,2,1^{n-5})$ where n is even.
  It is easy to see
that $\pi=(4,3^3,2,1)$  is potentially $(K_5-C_4)$-graphic. If
$n\geq 8$, then $(4,3^3,2,1^{n-5})$ has a realization containing a
$K_5-C_4$ (see Figure 11). Second, we consider
$\pi=(4,3^3,2^2,1^{n-6})$ where n is odd. It is easy to see that
$\pi=(4,3^3,2^2,1)$ is potentially $(K_5-C_4)$-graphic. If
$n\geq9$, then $(4,3^3,2^2,1^{n-6})$ has a realization containing
a $K_5-C_4$ (see Figure 12). Third, we consider
$\pi=(4,3^3,2^3,1^{n-7})$ where n is even. It is easy to see that
$\pi=(4,3^3,2^3,1)$ is potentially $(K_5-C_4)$-graphic. If
$n\geq10$, then $(4,3^3,2^3,1^{n-7})$ has a realization containing
a $K_5-C_4$ (see Figure 13). Then, we consider
$\pi=(4,3^3,2^b,1^{n-4-b})$ where $b\geq4$. In this case,
$(4,3^3,2^b,1^{n-4-b})$ has a realization containing a $K_5-C_4$
(see Figure 14). Thus, we are done.

\par
  Subcase 2.4: $a=4$. Then $\pi=(4,3^4,2^b,1^{n-5-b})$.
 There are two subcases: $b\geq1$ and $b=0$.
\par
  Suppose $b\geq1$. It is easy to see that
$\pi=(4,3^4,2,1^{n-6})$ and $\pi=(4,3^4,2^2,1^{n-7})$ are
potentially $(K_5-C_4)$-graphic (see Figure 15 and Figure 16,
respectively). If $b\geq3$, then $(4,3^4,2^b,1^{n-5-b})$ has a
realization containing a $K_5-C_4$ (see Figure 17).  Thus, we are
done.
\par
  Suppose $b=0$. Then $\pi=(4,3^4,1^{n-5})$. Since $\sigma(\pi)$ is even,
$n-5$ must be even. Clearly, $(4,3^4,1^{n-5})$ has a realization
containing a $K_5-C_4$ (see Figure 18). Thus, we are done.
\par
  Subcase 2.5:  $a\geq5$. Then  $\pi=(4,3^a,2^b,1^{n-1-a-b})$ where
$a\geq5$ and $n-1-a-b \geq 1$. There are two subcases: $b\geq1$
and $b=0$.
\par
  Suppose $b\geq1$.
\par
  If $a$ is even, it is easy to see that $\pi=(4,3^6,2,1^{n-8})$ has a realization containing a $K_5-C_4$ (see Figure 19)
. If $a=6$ and $b\geq2$, then $(4,3^6,2^b,1^{n-7-b})$ has a
realization containing a $K_5-C_4$ (see Figure 20). If $a\geq8$
and $b=1$, then $(4,3^a,2,1^{n-2-a})$ has a realization containing
a $K_5-C_4$ (see Figure 21). If $a\geq8$ and $b\geq2$, then
$(4,3^a,2^b,1^{n-1-a-b})$ has a realization containing a $K_5-C_4$
(see Figure 22).
\par
  If $a$ is odd, it is easy to see that
$\pi=(4,3^5,2,1^{n-7})$ has a realization containing a $K_5-C_4$
(see Figure 23). If $a=5$ and $b\geq2$, then
$(4,3^5,2^b,1^{n-6-b})$ has a realization containing a $K_5-C_4$
(see Figure 24). If $a\geq7$ and $b=1$, then $(4,3^a,2,1^{n-2-a})$
has a realization containing a $K_5-C_4$ (see Figure 25). If
$a\geq7$ and $b\geq2$, then $(4,3^a,2^b,1^{n-1-a-b})$ has a
realization containing a $K_5-C_4$ (see Figure 26). Thus, we are
done.
 \par
  Suppose $b=0$. Then $\pi=(4,3^a,1^{n-1-a})$. Since $\sigma(\pi)$ is even,
$n-1$ must be even. \par
  If $a$ is even, it is easy to see that
$\pi=(4,3^6,1^2)$ is potentially $(K_5-C_4)$-graphic. If $a=6$ and
$n\geq11$, then $(4,3^6,1^{n-7})$ has a realization containing a
$K_5-C_4$ (see Figure 27). If $a\geq8$, then $(4,3^a,1^{n-1-a})$
has a realization containing a $K_5-C_4$ (see Figure 28).
\par
  If $a$ is odd, it is easy to see that
$\pi=(4,3^5,1)$ and $\pi=(4,3^7,1)$ are potentially
$(K_5-C_4)$-graphic. If $a=5$ and $n\geq9$, then $(4,3^5,1^{n-6})$
has a realization containing a $K_5-C_4$ (see Figure 29). If $a=7$
and $n\geq11$, then $(4,3^7,1^{n-8})$ has a realization containing
a $K_5-C_4$ (see Figure 30). If $a\geq9$, then $(4,3^a,1^{n-1-a})$
has a realization containing a $K_5-C_4$ (see Figure 31). Thus, we
are done.\par

\par
  Subcase 3: $d_2=2$. Then $\pi=(4,2^a,1^{n-1-a})$ where $a\geq4$ and $n-1-a\geq1$. Since $\sigma(\pi)$ is even,
$n-1-a$ must be even. We are going to prove that $\pi$ is
potentially $(K_5-C_4)$-graphic. If $a=4$, then
$\pi=(4,2^4,1^{n-5})$ where $n-5$ is even. It is enough to show
$\pi_1=(1^{n-5})$ is graphic. Clearly, $\pi_1$ has a realization
consisting of $\frac {n-5}{2}$ disjoint edges. Thus,
$\pi=(4,2^4,1^{n-5})$ is potentially $(K_5-C_4)$-graphic. If
$a\geq5$, it is enough to show $\pi_1=(2^{a-4},1^{n-1-a})$ is
graphic. Clearly, $\pi_1$ is graphic. Thus,  we are done.\par
\par
\par
\par
\par
  \par
\section{  Application } \par

\par
  Using Theorem 3.1, we give a simple proof of the following theorem
due to  Lai [6]:
\par
  \textbf{Theorem 4.1 }  (Lai [6])  For $n\geq5$,
$\sigma(K_5-C_4,n)=4n-4$.
\par
Proof: First we claim that for $n\geq5,
\sigma(K_5-C_4,n)\geq4n-4$. It is enough to show that there exist
$\pi_1$ with $\sigma(\pi_1)=4n-6$, such that $\pi_1$ is not
potentially $(K_5-C_4)$-graphic. Take $\pi_1=((n-1)^2,2^{n-2})$,
then $\sigma(\pi_1)=4n-6$, and it is easy to see that $\pi_1$ is
not potentially $(K_5-C_4)$-graphic by Theorem 3.1.
\par
  Now we show that if $\pi$ is an $n$-term $(n\geq5)$ graphical
sequence with $\sigma(\pi)\geq4n-4$, then there exist a
realization of $\pi$ containing a $K_5-C_4$. Hence, it suffices to
show that $\pi$ is  potentially $(K_5-C_4)$-graphic.

\par
  If $d_5=1$, then $\sigma(\pi)=d_1+d_2+d_3+d_4+(n-4)$ and
$d_1+d_2+d_3+d_4\leq12+(n-4)=n+8$. Therefore,
$\sigma(\pi)\leq2n+4<4n-4$, which is a contradiction. Thus,
$d_5\geq2$.

\par
  If $d_1\leq3$, then $\sigma(\pi)\leq3n<4n-4$, which is a contradiction.
Thus, $d_1\geq4$.

\par
  Since $\sigma(\pi)\geq4n-4$, then $\pi$ is not one of the
following:
\par
  $((n-2)^2,2^{n-2})$ for $n\geq6$,
     $(n-k,k+i,2^i,1^{n-i-2})$ where $i=3,4,\cdots,n-2k$ and $k=1,2,\cdots,[\frac{n-1}{2}]-1$,
     $(4,2^5)$,
     $(4,2^6)$.
Thus, $\pi$ satisfies the conditions $(1)-(6)$ in Theorem 3.1.
Therefore, $\pi$ is potentially $(K_5-C_4)$-graphic.
\par

\par
\section * {Appendix}

\setlength{\unitlength}{1.2mm}
\begin{picture}(60,30)
\put(0, 20){\circle*{1}}
 \put(10, 20){\circle*{1}}
\put(20,20){\circle*{1}}
\put(0, 30){\circle*{1}}
 \put(10,30){\circle*{1}}
 \put(0, 30){\line(1, 0){10.0}}
\put(0, 20){\line(1, 0){10.0}}
\put(0, 20){\line(1, 1){10.0}}
\put(0, 30){\line(1, -1){10.0}}
\put(0, 20){\line(0, 1){10.0}}
\put(20, 20){\line(-1, 1){10.0}}
 \put(20, 20){\line(-2, 1){20.0}}
\put(10, 20){\line(0, 1){10.0}} \put(10, 20){\oval(20, 10)[b]}
\put(30, 20){\circle*{1}}
 \put(40, 20){\circle*{1}}
 \put(50, 20){\circle*{1}}
  \put(30,30){\circle*{1}}
  \put(40, 30){\circle*{1}}
   \put(50,30){\circle*{1}}
   \put(30, 30){\line(1, 0){10.0}}
   \put(30, 30){\line(0, -1){10.0}}
   \put(30, 20){\line(1, 0){10.0}}
   \put(40, 20){\line(1, 0){10.0}}
   \put(40, 30){\line(1, 0){10.0}}
   \put(40, 30){\line(0, -1){10.0}}
   \put(50, 30){\line(0, -1){10.0}}
   \multiput(51,20)(3,0){3}{\line(1,0){2}}
   \multiput(51,30)(3,0){3}{\line(1,0){2}}
   \put(60, 20){\circle*{1}}
  \put(60,30){\circle*{1}}
  \put(70, 20){\circle*{1}}
  \put(70,30){\circle*{1}}
  \put(60, 20){\line(1, 0){10.0}}
  \put(60, 30){\line(1, 0){10.0}}
  \put(60, 30){\line(0, -1){10.0}}
  \put(70, 30){\line(0, -1){10.0}}
  \put(30, 30){\line(4, -1){40.0}}
  \put(30, 20){\line(4, 1){40.0}}
  \put(25, 10){\makebox(8, 1)[l]{$(4^3, 3^{n-3})$}}
  \put(25, 5){\makebox(8, 1)[l]{Figure 1}}
  \put(0,18){$v_3$}
  \put(0,31){$v_2$}
  \put(10,18){$v_4$}
  \put(20,18){$v_5$}
  \put(10,31){$v_1$}
  \put(30,18){$v_7$}
  \put(30,31){$v_6$}
  \put(40,18){$v_9$}
  \put(50,18){$v_{11}$}
  \put(40,31){$v_8$}
  \put(60,18){$v_{n-2}$}
  \put(50,31){$v_{10}$}
  \put(70,18){$v_{n}$}
  \put(60,31){$v_{n-3}$}
  \put(70,31){$v_{n-1}$}

  \end{picture}
\vskip 0.1in \setlength{\unitlength}{1.2mm}
\begin{picture}(60,30)
\put(0, 20){\circle*{1}}
 \put(10, 20){\circle*{1}}
\put(20,20){\circle*{1}}
\put(0, 30){\circle*{1}}
 \put(10,30){\circle*{1}}
 \put(0, 30){\line(1, 0){10.0}}
\put(0, 20){\line(1, 0){10.0}}
\put(0, 20){\line(1, 1){10.0}}
\put(0, 30){\line(1, -1){10.0}}
\put(0, 20){\line(0, 1){10.0}}
\put(20, 20){\line(-1, 1){10.0}}
 \put(20, 20){\line(-2, 1){20.0}}
\put(10, 20){\line(0, 1){10.0}}
\put(20,30){\circle*{1}}
\put(30,30){\circle*{1}}
  \put(40, 30){\circle*{1}}
   \put(50,30){\circle*{1}}
   \put(30, 20){\circle*{1}}
 \put(40, 20){\circle*{1}}
 \put(50, 20){\circle*{1}}
 \put(20, 20){\line(0, 1){10.0}}
 \put(20, 30){\line(1, 0){10.0}}
 \put(30, 20){\line(0, 1){10.0}}
 \put(20, 30){\line(1, -1){10.0}}
 \put(30, 20){\line(1, 0){10.0}}
  \put(30, 30){\line(1, 0){10.0}}

    \put(40, 30){\line(0, -1){10.0}}
    \put(50, 30){\line(0, -1){10.0}}
    \multiput(41,20)(3,0){3}{\line(1,0){2}}
   \multiput(41,30)(3,0){3}{\line(1,0){2}}
   \put(60, 20){\circle*{1}}
  \put(60,30){\circle*{1}}
  \put(70, 20){\circle*{1}}
  \put(70,30){\circle*{1}}
  \put(60, 20){\line(1, 0){10.0}}
  \put(60, 30){\line(1, 0){10.0}}
  \put(50, 20){\line(1, 0){10.0}}
  \put(50, 30){\line(1, 0){10.0}}

  \put(70, 30){\line(0, -1){10.0}}
  \put(60, 20){\line(1, 1){10.0}}
  \put(60, 30){\line(1, -1){10.0}}
  \put(25, 10){\makebox(8, 1)[l]{$(4^2, 3^{n-2})$}}
  \put(25, 5){\makebox(8, 1)[l]{Figure 2}}
  \put(0,18){$v_3$}
  \put(0,31){$v_2$}
  \put(10,18){$v_4$}
  \put(20,18){$v_5$}
  \put(10,31){$v_1$}
  \put(20,31){$v_6$}
  \put(30,18){$v_8$}
  \put(30,31){$v_7$}
  \put(40,18){$v_{10}$}
  \put(50,18){$v_{n-4}$}
  \put(40,31){$v_9$}
  \put(60,18){$v_{n-2}$}
  \put(50,31){$v_{n-5}$}
  \put(70,18){$v_{n}$}
  \put(60,31){$v_{n-3}$}
  \put(70,31){$v_{n-1}$}
  \end{picture}
\vskip 0.1in \setlength{\unitlength}{1.2mm}
\begin{picture}(60,30)
\put(0, 20){\circle*{1}}
 \put(10, 20){\circle*{1}}
\put(20,20){\circle*{1}} \put(0, 30){\circle*{1}}
 \put(10,30){\circle*{1}}
 \put(0, 30){\line(1, 0){10.0}}
 \put(0, 20){\line(1, 1){10.0}}
\put(0, 30){\line(1, -1){10.0}} \put(0, 20){\line(0, 1){10.0}}
\put(20, 20){\line(-1, 1){10.0}} \put(10, 20){\line(1, 0){10.0}}

\put(10, 20){\line(0, 1){10.0}} \put(10, 20){\oval(20, 10)[b]}
\put(30, 20){\circle*{1}}
 \put(40, 20){\circle*{1}}
 \put(50, 20){\circle*{1}}
  \put(30,30){\circle*{1}}
  \put(40, 30){\circle*{1}}
   \put(50,30){\circle*{1}}
   \put(30, 30){\line(1, 0){10.0}}
   \put(30, 30){\line(0, -1){10.0}}
   \put(30, 20){\line(1, 0){10.0}}
   \put(40, 20){\line(1, 0){10.0}}
   \put(40, 30){\line(1, 0){10.0}}
   \put(40, 30){\line(0, -1){10.0}}
   \put(50, 30){\line(0, -1){10.0}}
   \multiput(51,20)(3,0){3}{\line(1,0){2}}
   \multiput(51,30)(3,0){3}{\line(1,0){2}}
   \put(60, 20){\circle*{1}}
  \put(60,30){\circle*{1}}
  \put(70, 20){\circle*{1}}
  \put(70,30){\circle*{1}}
  \put(60, 20){\line(1, 0){10.0}}
  \put(60, 30){\line(1, 0){10.0}}
  \put(60, 30){\line(0, -1){10.0}}
  \put(70, 30){\line(0, -1){10.0}}
  \put(30, 30){\line(4, -1){40.0}}
  \put(30, 20){\line(4, 1){40.0}}
  \put(25, 10){\makebox(8, 1)[l]{$(4, 3^{n-1})$}}
  \put(25, 5){\makebox(8, 1)[l]{Figure 3}}
  \put(0,18){$v_3$}
  \put(0,31){$v_2$}
  \put(10,18){$v_4$}
  \put(20,18){$v_5$}
  \put(10,31){$v_1$}
  \put(30,18){$v_7$}
  \put(30,31){$v_6$}
  \put(40,18){$v_9$}
  \put(50,18){$v_{11}$}
  \put(40,31){$v_8$}
  \put(60,18){$v_{n-2}$}
  \put(50,31){$v_{10}$}
  \put(70,18){$v_{n}$}
  \put(60,31){$v_{n-3}$}
  \put(70,31){$v_{n-1}$}
\end{picture}
\vskip 0.1in \setlength{\unitlength}{1.2mm}
\begin{picture}(60,30)
\put(0, 20){\circle*{1}}
 \put(10, 20){\circle*{1}}
\put(5,25){\circle*{1}}
 \put(0, 30){\circle*{1}}
 \put(10,30){\circle*{1}}
 \put(0, 20){\line(1, 0){10.0}}
 \put(0, 20){\line(1, 1){5.0}}
 \put(10, 20){\line(-1, 1){5.0}}
 \put(0, 30){\line(1, 0){10.0}}
 \put(0, 30){\line(1, -1){5.0}}
 \put(10, 30){\line(-1, -1){5.0}}
 \put(0, 20){\line(0, 1){10.0}}
 \put(20, 30){\circle*{1}}
 \put(30, 30){\circle*{1}}
 \put(40, 30){\circle*{1}}
 \put(50, 30){\circle*{1}}
\put(20,20){\circle*{1}} \put(50, 20){\circle*{1}}

 \put(10, 30){\line(1, 0){10.0}}

 \put(20, 20){\line(1, 0){30.0}}
 \put(20, 30){\line(1, 0){10.0}}
 \put(30, 30){\line(1, 0){10.0}}
 \put(10, 30){\line(1, -1){10.0}}
 \put(50, 30){\line(0, -1){10.0}}

 \multiput(41,30)(3,0){3}{\line(1,0){2}}

 \put(25, 10){\makebox(8, 1)[l]{$(4^2, 3^{2},2^{n-4})$}}
  \put(25, 5){\makebox(8, 1)[l]{Figure 4}}
  \put(0,31){$v_3$}
  \put(10,31){$v_2$}
  \put(0,18){$v_4$}
  \put(10,18){$v_5$}
  \put(5,23){$v_1$}
  \put(20,31){$v_7$}
  \put(20,18){$v_6$}
  \put(30,31){$v_8$}
  \put(40,31){$v_9$}
  \put(50,18){$v_{n}$}
  \put(50,31){$v_{n-1}$}
  \end{picture}
\vskip 0.1in \setlength{\unitlength}{1.2mm}
\begin{picture}(60,30)
\put(0, 20){\circle*{1}}
 \put(10, 20){\circle*{1}}
\put(5,25){\circle*{1}} \put(0, 30){\circle*{1}}
 \put(10,30){\circle*{1}}
 \put(0, 20){\line(1, 0){10.0}}
 \put(0, 20){\line(1, 1){5.0}}
 \put(10, 20){\line(-1, 1){5.0}}
 \put(0, 30){\line(1, 0){10.0}}
 \put(10, 30){\line(1, 0){10.0}}
 \put(0, 30){\line(1, -1){5.0}}
 \put(10, 30){\line(-1, -1){5.0}}
 \put(0, 20){\line(0, 1){10.0}}
 \put(10, 30){\line(0, -1){10.0}}
 \put(20, 20){\circle*{1}}
\put(20, 30){\circle*{1}}
 \put(30,30){\circle*{1}}
 \put(40, 30){\circle*{1}}
 \put(50,30){\circle*{1}}
 \put(50, 20){\circle*{1}}
\put(20, 30){\line(1, 0){10.0}} \put(30, 30){\line(1, 0){10.0}}
 \put(20, 20){\line(1, 0){30.0}}
\put(20, 30){\line(0, -1){10.0}} \put(50, 30){\line(0, -1){10.0}}
      \multiput(41,30)(3,0){3}{\line(1,0){2}}
\put(25, 10){\makebox(8, 1)[l]{$(4^2, 3^{4},2^{n-6})$}}
  \put(25, 5){\makebox(8, 1)[l]{Figure 5}}
  \end{picture}
\vskip 0.1in \setlength{\unitlength}{1.2mm}
\begin{picture}(60,30)
\put(0, 20){\circle*{1}}
 \put(10, 20){\circle*{1}}
\put(5,25){\circle*{1}} \put(0, 30){\circle*{1}}
 \put(10,30){\circle*{1}}
 \put(0, 20){\line(1, 0){10.0}}
 \put(0, 20){\line(1, 1){5.0}}
 \put(10, 20){\line(-1, 1){5.0}}
 \put(0, 30){\line(1, 0){10.0}}
 \put(0, 30){\line(1, -1){5.0}}
 \put(10, 30){\line(-1, -1){5.0}}
 \put(0, 20){\line(0, 1){10.0}}
 \put(10, 30){\line(0, -1){10.0}}
 \put(20, 20){\circle*{1}}
 \put(30, 20){\circle*{1}}
\put(20, 30){\circle*{1}}
 \put(30,30){\circle*{1}}
 \put(40, 30){\circle*{1}}
 \put(50,30){\circle*{1}}
 \put(40, 20){\circle*{1}}
 \put(60,30){\circle*{1}}
\put(70,30){\circle*{1}} \put(10,30){\line(1, 0){10.0}}
  \put(20, 30){\line(1, 0){10.0}}
   \put(30, 30){\line(1, 0){10.0}}
    \put(20, 20){\line(1, 0){10.0}}
\put(20, 30){\line(0, -1){10.0}}

      \put(30, 30){\line(-1, -1){10.0}}
      \put(40, 30){\line(-1, -1){10.0}}
      \put(50, 30){\line(-1, -1){10.0}}

      \put(50, 30){\line(1, 0){10.0}}
\put(70, 30){\line(-3, -1){30.0}}
\multiput(41,30)(3,0){3}{\line(1,0){2}}
      \multiput(31,20)(3,0){3}{\line(1,0){2}}
\multiput(61,30)(3,0){3}{\line(1,0){2}}
      \put(25, 10){\makebox(8, 1)[l]{$(4^2, 3^{a},2^{n-2-a})$}}
  \put(25, 5){\makebox(8, 1)[l]{Figure 6}}

  \end{picture}
\vskip 0.1in \setlength{\unitlength}{1.2mm}
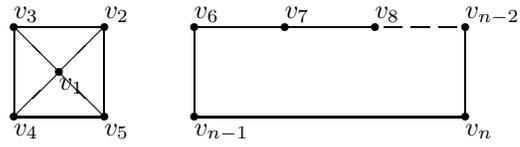
\begin{picture}(60,30)
\put(15, 20){\circle*{1}}
 \put(25, 20){\circle*{1}}
\put(20,25){\circle*{1}} \put(15, 30){\circle*{1}}
 \put(25,30){\circle*{1}}
 \put(15, 20){\line(1, 0){10.0}}
 \put(15, 20){\line(1, 1){5.0}}
 \put(25, 20){\line(-1, 1){5.0}}
 \put(15, 30){\line(1, 0){10.0}}
 \put(15, 30){\line(1, -1){5.0}}
 \put(25, 30){\line(-1, -1){5.0}}
 \put(15, 20){\line(0, 1){10.0}}
 \put(25, 30){\line(0, -1){10.0}}
 \put(35, 20){\circle*{1}}
\put(35, 30){\circle*{1}}
 \put(45,30){\circle*{1}}
 \put(55, 30){\circle*{1}}
 \put(65,30){\circle*{1}}
 \put(65, 20){\circle*{1}}
\put(35, 30){\line(1, 0){10.0}} \put(45, 30){\line(1, 0){10.0}}
 \put(35, 20){\line(1, 0){30.0}}
\put(35, 30){\line(0, -1){10.0}}
\put(65, 30){\line(0, -1){10.0}}
      \multiput(56,30)(3,0){3}{\line(1,0){2}}
\put(25, 10){\makebox(8, 1)[l]{$(4, 3^{4},2^{n-5})$}}
  \put(25, 5){\makebox(8, 1)[l]{Figure 7}}
  \put(15,31){$v_3$}
  \put(25,31){$v_2$}
  \put(15,18){$v_4$}
  \put(25,18){$v_5$}
  \put(20,23){$v_1$}
  \put(45,31){$v_7$}
  \put(35,31){$v_6$}
  \put(55,31){$v_8$}
  \put(65,18){$v_n$}
  \put(35,18){$v_{n-1}$}
  \put(65,31){$v_{n-2}$}
  \end{picture}
\vskip 0.1in \setlength{\unitlength}{1.2mm}
\begin{picture}(60,30)
\put(15, 20){\circle*{1}}
 \put(25, 20){\circle*{1}}
\put(20,25){\circle*{1}} \put(15, 30){\circle*{1}}
 \put(25,30){\circle*{1}}
 \put(15, 20){\line(1, 0){10.0}}
 \put(15, 20){\line(1, 1){5.0}}
 \put(25, 20){\line(-1, 1){5.0}}
 \put(15, 30){\line(1, 0){10.0}}

 \put(15, 30){\line(1, -1){5.0}}
 \put(25, 30){\line(-1, -1){5.0}}
 \put(15, 20){\line(0, 1){10.0}}
 \put(25, 30){\line(0, -1){10.0}}
 \put(35, 20){\circle*{1}}
\put(35, 30){\circle*{1}}
 \put(45,30){\circle*{1}}
 \put(55, 30){\circle*{1}}
 \put(65,30){\circle*{1}}
 \put(45, 20){\circle*{1}}
\put(35, 30){\line(1, 0){10.0}} \put(45, 30){\line(1, 0){10.0}}
\put(35, 30){\line(1, -1){10.0}}
 \put(35, 20){\line(1, 0){10.0}}
\put(35, 30){\line(0, -1){10.0}} \put(65, 30){\line(-2, -1){20.0}}
      \multiput(56,30)(3,0){3}{\line(1,0){2}}
\put(25, 10){\makebox(8, 1)[l]{$(4, 3^{6},2^{n-7})$}}
  \put(25, 5){\makebox(8, 1)[l]{Figure 8}}
  \end{picture}
\vskip 0.1in \setlength{\unitlength}{1.2mm}
\begin{picture}(60,30)
\put(5, 20){\circle*{1}}
 \put(15, 20){\circle*{1}}
\put(10,25){\circle*{1}} \put(5, 30){\circle*{1}}
 \put(15,30){\circle*{1}}
 \put(5, 20){\line(1, 0){10.0}}
 \put(5, 20){\line(1, 1){5.0}}
 \put(15, 20){\line(-1, 1){5.0}}
 \put(5, 30){\line(1, 0){10.0}}
 \put(5, 30){\line(1, -1){5.0}}
 \put(15, 30){\line(-1, -1){5.0}}
 \put(5, 20){\line(0, 1){10.0}}
 \put(15, 30){\line(0, -1){10.0}}
 \put(25, 20){\circle*{1}}
 \put(35, 20){\circle*{1}}
\put(45,20){\circle*{1}}
 \put(25, 30){\circle*{1}}
 \put(35,30){\circle*{1}}
 \put(45, 30){\circle*{1}}
 \put(55,30){\circle*{1}}
 \put(55,20){\circle*{1}}
\put(65,30){\circle*{1}} \put(55, 30){\line(1, 0){10.0}} \put(25,
30){\line(1, 0){10.0}}
   \put(35, 30){\line(1, -1){10.0}}
   \put(45, 30){\line(1, -1){10.0}}
   \put(25, 30){\line(1, -1){10.0}}
   \put(55, 30){\line(0, -1){10.0}}
\put(25, 20){\line(1, 0){10.0}} \put(45, 30){\line(1, 0){10.0}}
\put(25, 20){\line(4, 1){40.0}} \put(35, 20){\line(1, 0){10.0}}
\put(25, 30){\line(0, -1){10.0}}
\multiput(36,30)(3,0){3}{\line(1,0){2}}
\multiput(46,20)(3,0){3}{\line(1,0){2}} \put(25,
10){\makebox(8,1)[l]{$(4, 3^{a},2)$}}
  \put(25, 5){\makebox(8, 1)[l]{Figure 9}}
  \end{picture}
\vskip 0.1in \setlength{\unitlength}{1.2mm}
\begin{picture}(60,30)
\put(5, 20){\circle*{1}}
 \put(15, 20){\circle*{1}}
\put(10,25){\circle*{1}} \put(5, 30){\circle*{1}}
 \put(15,30){\circle*{1}}
 \put(5, 20){\line(1, 0){10.0}}
 \put(5, 20){\line(1, 1){5.0}}
 \put(15, 20){\line(-1, 1){5.0}}
 \put(5, 30){\line(1, 0){10.0}}
 \put(5, 30){\line(1, -1){5.0}}
 \put(15, 30){\line(-1, -1){5.0}}
 \put(5, 20){\line(0, 1){10.0}}
 \put(15, 30){\line(0, -1){10.0}}
 \put(25, 20){\circle*{1}}
 \put(35, 20){\circle*{1}}
\put(45,20){\circle*{1}}
 \put(25, 30){\circle*{1}}
 \put(35,30){\circle*{1}}
 \put(45, 30){\circle*{1}}
 \put(55,30){\circle*{1}}
 \put(55,20){\circle*{1}}

 \put(65,30){\circle*{1}}
 \put(75, 30){\circle*{1}}
\put(55, 30){\line(1, 0){10.0}}

   \put(25, 30){\line(1, 0){10.0}}
   \put(35, 30){\line(1, -1){10.0}}
   \put(45, 30){\line(1, -1){10.0}}
   \put(25, 30){\line(1, -1){10.0}}
   \put(55, 30){\line(0, -1){10.0}}
\put(25, 20){\line(1, 0){10.0}} \put(45, 30){\line(1, 0){10.0}}

     \put(25, 20){\line(5, 1){50.0}}
\put(35, 20){\line(1, 0){10.0}}
      \put(25, 30){\line(0, -1){10.0}}
\multiput(36,30)(3,0){3}{\line(1,0){2}}
      \multiput(66,30)(3,0){3}{\line(1,0){2}}
      \multiput(46,20)(3,0){3}{\line(1,0){2}}
\put(25, 10){\makebox(8, 1)[l]{$(4, 3^{a},2^{n-1-a})$}}
  \put(25, 5){\makebox(8, 1)[l]{Figure 10}}

  \end{picture}
\vskip 0.1in \setlength{\unitlength}{1.2mm}
\begin{picture}(60,30)
\put(5, 20){\circle*{1}}
 \put(15, 20){\circle*{1}}
\put(10,25){\circle*{1}} \put(5, 30){\circle*{1}}
 \put(15,30){\circle*{1}}
 \put(5, 20){\line(1, 0){10.0}}
 \put(5, 20){\line(1, 1){5.0}}
 \put(15, 20){\line(-1, 1){5.0}}
 \put(5, 30){\line(1, 0){10.0}}
 \put(5, 30){\line(1, -1){5.0}}
 \put(15, 30){\line(-1, -1){5.0}}
 \put(5, 20){\line(0, 1){10.0}}
 \put(15, 30){\line(1, 0){10.0}}

 \put(35, 20){\circle*{1}}
\put(45,20){\circle*{1}}
 \put(25, 30){\circle*{1}}
 \put(35,30){\circle*{1}}
 \put(45, 30){\circle*{1}}
 \put(55,30){\circle*{1}}
 \put(55, 20){\circle*{1}}
 \put(65,30){\circle*{1}}
 \put(65, 20){\circle*{1}}
 \put(45, 30){\line(0, -1){10.0}}
 \put(55, 30){\line(0, -1){10.0}}

\put(35, 30){\line(0, -1){10.0}} \put(65, 30){\line(0, -1){10.0}}
\multiput(58,25)(2,0){3}{\circle*{0.3}{6}}

\put(25, 10){\makebox(8, 1)[l]{$(4, 3^{3},2,1^{n-5})$}}
\put(25,5){\makebox(8, 1)[l]{Figure 11}} \put(5,31){$v_3$}
  \put(15,31){$v_2$}
  \put(5,18){$v_4$}
  \put(15,18){$v_5$}
  \put(10,23){$v_1$}
  \put(35,31){$v_7$}
  \put(25,31){$v_6$}
  \put(55,18){$v_{12}$}
  \put(45,18){$v_{10}$}
\put(45,31){$v_9$}
  \put(65,18){$v_n$}
  \put(55,18){$v_{12}$}
  \put(55,31){$v_{11}$}
  \put(65,31){$v_{n-1}$}
\end{picture}
\vskip 0.1in \setlength{\unitlength}{1.2mm}
\begin{picture}(60,30)
\put(0, 20){\circle*{1}}
 \put(10, 20){\circle*{1}}
\put(5,25){\circle*{1}} \put(0, 30){\circle*{1}}
 \put(10,30){\circle*{1}}
 \put(0, 20){\line(1, 0){10.0}}
 \put(0, 20){\line(1, 1){5.0}}
 \put(10, 20){\line(-1, 1){5.0}}
 \put(10, 30){\line(1, 0){5.0}}
 \put(10, 30){\line(-1, 0){10.0}}
 \put(0, 30){\line(1, -1){5.0}}
 \put(10, 30){\line(-1, -1){5.0}}
 \put(0, 20){\line(0, 1){10.0}}
\put(20, 30){\circle*{1}} \put(30,20){\circle*{1}}
 \put(20, 20){\circle*{1}}
 \put(30,30){\circle*{1}}
\put(40, 20){\circle*{1}}
\put(40, 30){\circle*{1}}
 \put(50,30){\circle*{1}}
 \put(50, 20){\circle*{1}}
\put(20, 30){\line(0, -1){10.0}}
 \put(30, 30){\line(0, -1){10.0}}
 \put(40, 30){\line(0, -1){10.0}}
 \put(50, 30){\line(0, -1){10.0}}
 \put(10, 30){\line(1, 0){10.0}}
 \multiput(43,25)(2,0){3}{\circle*{0.3}{3}}
\put(20, 10){\makebox(8, 1)[l]{$(4, 3^{3},2^2,1^{n-6})$}}
\put(20,5){\makebox(8, 1)[l]{Figure 12}}
 \put(0,31){$v_3$}
  \put(10,31){$v_2$}
  \put(0,18){$v_4$}
  \put(10,18){$v_5$}
  \put(5,23){$v_1$}
  \put(20,18){$v_7$}
  \put(20,31){$v_6$}
  \put(30,31){$v_8$}
  \put(40,31){$v_{10}$}
\put(30,18){$v_9$}
  \put(50,18){$v_n$}
  \put(40,18){$v_{11}$}
  \put(50,31){$v_{n-1}$}
\end{picture}
\vskip 0.1in \setlength{\unitlength}{1.2mm}
\begin{picture}(60,30)
\put(0, 20){\circle*{1}}
 \put(10, 20){\circle*{1}}
\put(5,25){\circle*{1}} \put(0, 30){\circle*{1}}
 \put(10,30){\circle*{1}}
 \put(0, 20){\line(1, 0){10.0}}
 \put(0, 20){\line(1, 1){5.0}}
 \put(10, 20){\line(-1, 1){5.0}}
 \put(0, 30){\line(1, 0){10.0}}
 \put(0, 30){\line(1, -1){5.0}}
 \put(10, 30){\line(-1, -1){5.0}}
 \put(0, 20){\line(0, 1){10.0}}
 \put(10, 30){\line(1, 0){10.0}}
\put(30,20){\circle*{1}}
 \put(40, 20){\circle*{1}}
 \put(30,30){\circle*{1}}
 \put(20,30){\circle*{1}}
 \put(40, 30){\circle*{1}}
 \put(50,30){\circle*{1}}
 \put(50, 20){\circle*{1}}
 \put(60,30){\circle*{1}}
 \put(60, 20){\circle*{1}}
 \put(40, 30){\line(0, -1){10.0}}
 \put(30, 30){\line(0, -1){10.0}}
\put(20, 30){\line(1, 0){10.0}}
 \put(50, 30){\line(0, -1){10.0}}
\put(60, 30){\line(0, -1){10.0}}
 \multiput(53,25)(2,0){3}{\circle*{0.3}{3}}
\put(20, 10){\makebox(8, 1)[l]{$(4, 3^{3},2^3,1^{n-7})$}}
\put(20,5){\makebox(8, 1)[l]{Figure 13}} \put(0,31){$v_3$}
  \put(10,31){$v_2$}
  \put(0,18){$v_4$}
  \put(10,18){$v_5$}
  \put(5,23){$v_1$}
\put(20,31){$v_6$}
  \put(30,31){$v_7$}
  \put(40,31){$v_{9}$}
  \put(50,31){$v_{11}$}
\put(30,18){$v_8$}
  \put(60,18){$v_n$}
  \put(40,18){$v_{10}$}
  \put(50,18){$v_{12}$}
  \put(60,31){$v_{n-1}$}
\end{picture}
\vskip 0.1in \setlength{\unitlength}{1.2mm}
\begin{picture}(60,30)
\put(0, 20){\circle*{1}}
 \put(10, 20){\circle*{1}}
\put(5,25){\circle*{1}} \put(0, 30){\circle*{1}}
 \put(10,30){\circle*{1}}
 \put(0, 20){\line(1, 0){10.0}}
 \put(0, 20){\line(1, 1){5.0}}
 \put(10, 20){\line(-1, 1){5.0}}
\put(10, 30){\line(-1, 0){10.0}}
\put(0, 30){\line(1, -1){5.0}}
 \put(10, 30){\line(-1, -1){5.0}}
 \put(0, 20){\line(0, 1){10.0}}
\put(30,20){\circle*{1}}
\put(20,30){\circle*{1}}
 \put(30,30){\circle*{1}}
 \put(40, 30){\circle*{1}}
 \put(50,30){\circle*{1}}
 \put(60, 20){\circle*{1}}
\put(30, 30){\line(0, -1){10.0}}
 \put(60, 30){\line(0, -1){10.0}}
 \put(30, 30){\line(1, 0){10.0}}
 \put(10, 30){\line(1, 0){10.0}}
 \put(40, 30){\line(1, 0){10.0}}
\put(30, 20){\line(1, 0){30.0}}
 \put(70, 30){\line(0, -1){10.0}}
\put(70,30){\circle*{1}}
 \put(80, 20){\circle*{1}}
 \put(80,30){\circle*{1}}
 \put(60,30){\circle*{1}}
 \put(70,20){\circle*{1}}
\put(80, 30){\line(0, -1){10.0}}
\multiput(51,30)(3,0){3}{\line(1,0){2}}
\multiput(73,25)(2,0){3}{\circle*{0.3}{3}}

\put(20,10){\makebox(8,1)[l]{$(4, 3^{3},2^b,1^{n-4-b})$}}
\put(20,5){\makebox(8, 1)[l]{Figure 14}}
\end{picture}
\vskip 0.1in \setlength{\unitlength}{1.2mm}
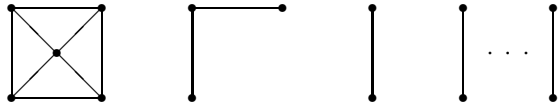
\begin{picture}(60,30)
\put(0, 20){\circle*{1}}
 \put(10, 20){\circle*{1}}
\put(5,25){\circle*{1}} \put(0, 30){\circle*{1}}
 \put(10,30){\circle*{1}}
 \put(0, 20){\line(1, 0){10.0}}
 \put(0, 20){\line(1, 1){5.0}}
 \put(10, 20){\line(-1, 1){5.0}}
 \put(0, 30){\line(1, 0){10.0}}
 \put(0, 30){\line(1, -1){5.0}}
 \put(10, 30){\line(-1, -1){5.0}}
 \put(0, 20){\line(0, 1){10.0}}
 \put(10, 30){\line(0, -1){10.0}}
\put(20,20){\circle*{1}} \put(20,30){\circle*{1}}
 \put(40, 20){\circle*{1}}
 \put(30,30){\circle*{1}}
 \put(40, 30){\circle*{1}}
 \put(50,30){\circle*{1}}
 \put(50, 20){\circle*{1}}
 \put(60,30){\circle*{1}}
 \put(60, 20){\circle*{1}}
 \put(40, 30){\line(0, -1){10.0}}
 \put(50, 30){\line(0, -1){10.0}}
\put(60, 30){\line(0, -1){10.0}}
 \put(20, 30){\line(0, -1){10.0}}
\put(20, 30){\line(1, 0){10.0}}
\multiput(53,25)(2,0){3}{\circle*{0.3}{3}}

\put(20,10){\makebox(8,1)[l]{$(4, 3^{4},2,1^{n-6})$}}
\put(20,5){\makebox(8, 1)[l]{Figure 15}}
\end{picture}
\vskip 0.1in \setlength{\unitlength}{1.2mm}
\begin{picture}(60,30)
\put(0, 20){\circle*{1}}
 \put(10, 20){\circle*{1}}
\put(5,25){\circle*{1}} \put(0, 30){\circle*{1}}
 \put(10,30){\circle*{1}}
 \put(0, 20){\line(1, 0){10.0}}
 \put(0, 20){\line(1, 1){5.0}}
 \put(10, 20){\line(-1, 1){5.0}}
 \put(0, 30){\line(1, 0){10.0}}
 \put(0, 30){\line(1, -1){5.0}}
 \put(10, 30){\line(-1, -1){5.0}}
 \put(0, 20){\line(0, 1){10.0}}
 \put(10, 30){\line(0, -1){10.0}}
\put(20,20){\circle*{1}}
 \put(20,30){\circle*{1}}
 \put(30,20){\circle*{1}}
 \put(40, 20){\circle*{1}}
 \put(30,30){\circle*{1}}
 \put(40, 30){\circle*{1}}
 \put(50,30){\circle*{1}}
 \put(50, 20){\circle*{1}}
 \put(60,30){\circle*{1}}
 \put(60, 20){\circle*{1}}
 \put(40, 30){\line(0, -1){10.0}}
 \put(50, 30){\line(0, -1){10.0}}
 \put(30, 30){\line(0, -1){10.0}}
\put(60, 30){\line(0, -1){10.0}}
 \put(20, 30){\line(0, -1){10.0}}
\put(20, 30){\line(1, 0){10.0}}
\multiput(53,25)(2,0){3}{\circle*{0.3}{3}}
\put(20,10){\makebox(8,1)[l]{$(4, 3^{4},2^2,1^{n-7})$}}
\put(20,5){\makebox(8, 1)[l]{Figure 16}}
\end{picture}
\vskip 0.1in \setlength{\unitlength}{1.2mm}
\begin{picture}(60,30)
\put(0, 20){\circle*{1}}
 \put(10, 20){\circle*{1}}
\put(5,25){\circle*{1}} \put(0, 30){\circle*{1}}
 \put(10,30){\circle*{1}}
 \put(0, 20){\line(1, 0){10.0}}
 \put(0, 20){\line(1, 1){5.0}}
 \put(10, 20){\line(-1, 1){5.0}}
 \put(0, 30){\line(1, 0){10.0}}
 \put(0, 30){\line(1, -1){5.0}}
 \put(10, 30){\line(-1, -1){5.0}}
 \put(0, 20){\line(0, 1){10.0}}
 \put(10, 30){\line(0, -1){10.0}}
\put(20,20){\circle*{1}} \put(20,30){\circle*{1}}
\put(30,30){\circle*{1}} \put(40,30){\circle*{1}}
\put(50,30){\circle*{1}} \put(50, 20){\circle*{1}}
 \put(60,30){\circle*{1}}
 \put(60, 20){\circle*{1}}
\put(20, 30){\line(0, -1){10.0}}
 \put(50, 30){\line(0, -1){10.0}}
\put(20, 30){\line(1, 0){10.0}} \put(30, 30){\line(1, 0){10.0}}
\put(20, 20){\line(1, 0){30.0}} \put(60, 30){\line(0, -1){10.0}}
\put(70,30){\circle*{1}}
 \put(70, 20){\circle*{1}}
\put(80,30){\circle*{1}}
 \put(80, 20){\circle*{1}}
 \put(70, 30){\line(0, -1){10.0}}
\put(80, 30){\line(0, -1){10.0}}
\multiput(41,30)(3,0){3}{\line(1,0){2}}
\multiput(73,25)(2,0){3}{\circle*{0.3}{3}}

\put(20,10){\makebox(8,1)[l]{$(4, 3^{4},2^b,1^{n-5-b})$}}
\put(20,5){\makebox(8, 1)[l]{Figure 17}}
\end{picture}
\vskip 0.1in \setlength{\unitlength}{1.2mm}
\begin{picture}(60,30)
\put(0, 20){\circle*{1}}
 \put(10, 20){\circle*{1}}
\put(5,25){\circle*{1}} \put(0, 30){\circle*{1}}
 \put(10,30){\circle*{1}}
 \put(0, 20){\line(1, 0){10.0}}
 \put(0, 20){\line(1, 1){5.0}}
 \put(10, 20){\line(-1, 1){5.0}}
 \put(0, 30){\line(1, 0){10.0}}
 \put(0, 30){\line(1, -1){5.0}}
 \put(10, 30){\line(-1, -1){5.0}}
 \put(0, 20){\line(0, 1){10.0}}
 \put(10, 30){\line(0, -1){10.0}}
\put(20,20){\circle*{1}}
 \put(20,30){\circle*{1}}
 \put(30,20){\circle*{1}}
 \put(40, 20){\circle*{1}}
 \put(30,30){\circle*{1}}
 \put(40, 30){\circle*{1}}
 \put(30, 30){\line(0, -1){10.0}}
 \put(40, 30){\line(0, -1){10.0}}
\put(20, 30){\line(0, -1){10.0}}
\multiput(33,25)(2,0){3}{\circle*{0.3}{3}}
\put(20,10){\makebox(8,1)[l]{$(4, 3^{4},1^{n-5})$}}
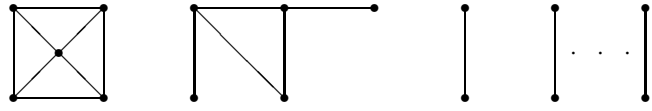
\put(20,5){\makebox(8, 1)[l]{Figure 18}}
\end{picture}
\vskip 0.1in \setlength{\unitlength}{1.2mm}
\begin{picture}(60,30)
\put(0, 20){\circle*{1}}
 \put(10, 20){\circle*{1}}
\put(5,25){\circle*{1}} \put(0, 30){\circle*{1}}
 \put(10,30){\circle*{1}}
 \put(50, 20){\circle*{1}}
 \put(0, 20){\line(1, 0){10.0}}
 \put(0, 20){\line(1, 1){5.0}}
 \put(10, 20){\line(-1, 1){5.0}}
 \put(0, 30){\line(1, 0){10.0}}
\put(0, 30){\line(1, -1){5.0}}
 \put(10, 30){\line(-1, -1){5.0}}
 \put(0, 20){\line(0, 1){10.0}}
 \put(10, 30){\line(0, -1){10.0}}
 \put(20, 20){\circle*{1}}
\put(20, 30){\circle*{1}}
 \put(30,30){\circle*{1}}
 \put(30, 20){\circle*{1}}
 \put(40, 30){\circle*{1}}
 \put(50,30){\circle*{1}}
 \put(70, 20){\circle*{1}}
 \put(60, 20){\circle*{1}}
\put(60, 30){\circle*{1}}
 \put(70,30){\circle*{1}}
 \put(20, 30){\line(1, 0){10.0}}
  \put(30,30){\line(1, 0){10.0}}
  \put(20, 30){\line(1, -1){10.0}}
   \put(30,30){\line(0, -1){10.0}}
 \put(50,30){\line(0, -1){10.0}}
\put(20, 30){\line(0, -1){10.0}} \put(60, 30){\line(0, -1){10.0}}
\put(70, 30){\line(0, -1){10.0}}
\multiput(62,25)(3,0){3}{\circle*{0.3}{2}} \put(25,
10){\makebox(8, 1)[l]{$(4, 3^{6},2,1^{n-8})$}}
  \put(25, 5){\makebox(8, 1)[l]{Figure 19}}
  \end{picture}
\vskip 0.1in \setlength{\unitlength}{1.2mm}
\begin{picture}(60,30)
\put(0, 20){\circle*{1}}
 \put(10, 20){\circle*{1}}
\put(5,25){\circle*{1}} \put(0, 30){\circle*{1}}
 \put(10,30){\circle*{1}}
 \put(0, 20){\line(1, 0){10.0}}
 \put(0, 20){\line(1, 1){5.0}}
 \put(10, 20){\line(-1, 1){5.0}}
 \put(0, 30){\line(1, 0){10.0}}
\put(0, 30){\line(1, -1){5.0}}
 \put(10, 30){\line(-1, -1){5.0}}
 \put(0, 20){\line(0, 1){10.0}}
 \put(10, 30){\line(0, -1){10.0}}
 \put(20, 20){\circle*{1}}
\put(20, 30){\circle*{1}}
 \put(30,30){\circle*{1}}
 \put(30, 20){\circle*{1}}
 \put(40, 30){\circle*{1}}
 \put(50,30){\circle*{1}}
 \put(70, 20){\circle*{1}}
 \put(60, 20){\circle*{1}}
\put(60, 30){\circle*{1}}
 \put(70,30){\circle*{1}}
 \put(80, 30){\circle*{1}}
\put(80, 20){\circle*{1}} \put(20, 30){\line(1, 0){10.0}} \put(30,
30){\line(1, 0){10.0}} \put(20, 30){\line(1, -1){10.0}}
 \put(20, 20){\line(1, 0){10.0}}
\put(20, 30){\line(0, -1){10.0}} \put(60, 30){\line(0, -1){10.0}}
\put(70, 30){\line(0, -1){10.0}} \put(80, 30){\line(0, -1){10.0}}
 \put(50, 30){\line(-2, -1){20.0}}
      \multiput(41,30)(3,0){3}{\line(1,0){2}}
      \multiput(72,25)(3,0){3}{\circle*{0.3}{2}}
\put(25, 10){\makebox(8, 1)[l]{$(4, 3^{6},2^b,1^{n-7-b})$}}
  \put(25, 5){\makebox(8, 1)[l]{Figure 20}}
  \end{picture}
\vskip 0.1in \setlength{\unitlength}{1.2mm}
\begin{picture}(60,30)
\put(0, 20){\circle*{1}}
 \put(10, 20){\circle*{1}}
\put(5,25){\circle*{1}} \put(0, 30){\circle*{1}}
 \put(10,30){\circle*{1}}
 \put(0, 20){\line(1, 0){10.0}}
 \put(0, 20){\line(1, 1){5.0}}
 \put(10, 20){\line(-1, 1){5.0}}
 \put(0, 30){\line(1, 0){10.0}}
 \put(0, 30){\line(1, -1){5.0}}
 \put(10, 30){\line(-1, -1){5.0}}
 \put(0, 20){\line(0, 1){10.0}}
 \put(10, 30){\line(0, -1){10.0}}
 \put(20, 20){\circle*{1}}
 \put(30, 20){\circle*{1}}
\put(40,20){\circle*{1}}
 \put(20, 30){\circle*{1}}
 \put(30,30){\circle*{1}}
 \put(40, 30){\circle*{1}}
 \put(50,30){\circle*{1}}
\put(50, 20){\circle*{1}}
 \put(60,30){\circle*{1}}
 \put(70, 30){\circle*{1}}
\put(80,30){\circle*{1}} \put(80, 20){\circle*{1}}
\put(70,20){\circle*{1}} \put(80, 30){\line(0, -1){10.0}}
 \put(50, 30){\line(0, -1){10.0}}
 \put(70, 30){\line(0, -1){10.0}}
\multiput(73,25)(2,0){3}{\circle*{0.3}{3}}
\put(20,30){\line(1,0){10.0}} \put(30, 30){\line(1, -1){10.0}}
\put(40, 30){\line(1, -1){10.0}}
   \put(20, 30){\line(1, -1){10.0}}
    \put(20, 20){\line(1, 0){10.0}}
    \put(50, 30){\line(1, 0){10.0}}
    \put(20, 20){\line(4, 1){40.0}}
\put(40, 30){\line(1, 0){10.0}} \put(30, 20){\line(1, 0){10.0}}
      \put(20, 30){\line(0, -1){10.0}}

\multiput(31,30)(3,0){3}{\line(1,0){2}}
\multiput(41,20)(3,0){3}{\line(1,0){2}}
\put(20,10){\makebox(8,1)[l]{$(4, 3^{a},2,1^{n-2-a})$}}
\put(20,5){\makebox(8, 1)[l]{Figure 21}}
\end{picture}
\vskip 0.1in \setlength{\unitlength}{1.2mm}
\begin{picture}(60,30)
\put(0, 20){\circle*{1}}
 \put(10, 20){\circle*{1}}
\put(5,25){\circle*{1}} \put(0, 30){\circle*{1}}
 \put(10,30){\circle*{1}}
 \put(0, 20){\line(1, 0){10.0}}
 \put(0, 20){\line(1, 1){5.0}}
 \put(10, 20){\line(-1, 1){5.0}}
 \put(0, 30){\line(1, 0){10.0}}
 \put(0, 30){\line(1, -1){5.0}}
 \put(10, 30){\line(-1, -1){5.0}}
 \put(0, 20){\line(0, 1){10.0}}
 \put(10, 30){\line(0, -1){10.0}}
 \put(20, 20){\circle*{1}}
 \put(30, 20){\circle*{1}}
\put(40,20){\circle*{1}}
 \put(20, 30){\circle*{1}}
 \put(30,30){\circle*{1}}
 \put(40, 30){\circle*{1}}
 \put(50,30){\circle*{1}}
\put(50, 20){\circle*{1}}
 \put(60,30){\circle*{1}}
 \put(70, 30){\circle*{1}}
\put(80,30){\circle*{1}}
 \put(90, 30){\circle*{1}}
 \put(80, 20){\circle*{1}}
 \put(90,20){\circle*{1}}
\put(80, 30){\line(0, -1){10.0}} \put(50, 30){\line(0, -1){10.0}}
 \put(90, 30){\line(0, -1){10.0}}
\multiput(83,25)(2,0){3}{\circle*{0.3}{3}}
\put(20,30){\line(1,0){10.0}} \put(30, 30){\line(1, -1){10.0}}
\put(40, 30){\line(1, -1){10.0}}
   \put(20, 30){\line(1, -1){10.0}}
    \put(20, 20){\line(1, 0){10.0}}
    \put(50, 30){\line(1, 0){10.0}}
    \put(20, 20){\line(5, 1){50.0}}
\put(40, 30){\line(1, 0){10.0}} \put(30, 20){\line(1, 0){10.0}}
      \put(20, 30){\line(0, -1){10.0}}

\multiput(31,30)(3,0){3}{\line(1,0){2}}
      \multiput(61,30)(3,0){3}{\line(1,0){2}}
      \multiput(41,20)(3,0){3}{\line(1,0){2}}
\put(20,10){\makebox(8,1)[l]{$(4, 3^{a},2^b,1^{n-1-a-b})$}}
\put(20,5){\makebox(8, 1)[l]{Figure 22}}
\end{picture}
\vskip 0.1in \setlength{\unitlength}{1.2mm}
\begin{picture}(60,30)
\put(0, 20){\circle*{1}}
 \put(10, 20){\circle*{1}}
\put(5,25){\circle*{1}} \put(0, 30){\circle*{1}}
 \put(10,30){\circle*{1}}
 \put(0, 20){\line(1, 0){10.0}}
 \put(0, 20){\line(1, 1){5.0}}
 \put(10, 20){\line(-1, 1){5.0}}
 \put(0, 30){\line(1, 0){10.0}}
 \put(0, 30){\line(1, -1){5.0}}
 \put(10, 30){\line(-1, -1){5.0}}
 \put(0, 20){\line(0, 1){10.0}}
 \put(50, 30){\line(0, -1){10.0}}
  \put(10, 20){\line(1, 0){10.0}}
   \put(10, 30){\line(1, 0){10.0}}
\put(20,20){\circle*{1}} \put(20,30){\circle*{1}}
\put(30,30){\circle*{1}} \put(40,30){\circle*{1}}
\put(50,30){\circle*{1}}
 \put(60,30){\circle*{1}}
\put(20, 30){\line(0, -1){10.0}} \put(20, 30){\line(1, 0){10.0}}
\put(40, 30){\line(0, -1){10.0}} \put(50,20){\circle*{1}}
 \put(40, 20){\circle*{1}}
 \put(60,20){\circle*{1}}
 \put(60, 30){\line(0, -1){10.0}}
\multiput(53,25)(2,0){3}{\circle*{0.3}{3}}

\put(20,10){\makebox(8,1)[l]{$(4, 3^{5},2,1^{n-7})$}}
\put(20,5){\makebox(8, 1)[l]{Figure 23}}
\end{picture}
\vskip 0.1in \setlength{\unitlength}{1.2mm}
\begin{picture}(60,30)
\put(0, 20){\circle*{1}}
 \put(10, 20){\circle*{1}}
\put(5,25){\circle*{1}} \put(0, 30){\circle*{1}}
 \put(10,30){\circle*{1}}
 \put(0, 20){\line(1, 0){10.0}}
 \put(0, 20){\line(1, 1){5.0}}
 \put(10, 20){\line(-1, 1){5.0}}
 \put(0, 30){\line(1, 0){10.0}}
 \put(0, 30){\line(1, -1){5.0}}
 \put(10, 30){\line(-1, -1){5.0}}
 \put(0, 20){\line(0, 1){10.0}}
 \put(10, 30){\line(0, -1){10.0}}
\put(20,20){\circle*{1}} \put(20,30){\circle*{1}}
\put(30,30){\circle*{1}} \put(40,30){\circle*{1}}
\put(50,30){\circle*{1}}
 \put(60,30){\circle*{1}}
\put(20, 30){\line(0, -1){10.0}} \put(20, 30){\line(1, 0){10.0}}
 \put(30, 30){\line(1, 0){10.0}}
 \put(50, 30){\line(1, 0){10.0}}
\put(20, 20){\line(1, 1){10.0}} \put(70,30){\circle*{1}}
 \put(70, 20){\circle*{1}}
\put(80,30){\circle*{1}}
 \put(80, 20){\circle*{1}}
 \put(70, 30){\line(0, -1){10.0}}
\put(80, 30){\line(0, -1){10.0}}
\multiput(41,30)(3,0){3}{\line(1,0){2}}
\multiput(73,25)(2,0){3}{\circle*{0.3}{3}}

\put(20,10){\makebox(8,1)[l]{$(4, 3^{5},2^b,1^{n-6-b})$}}
\put(20,5){\makebox(8, 1)[l]{Figure 24}}
\end{picture}
\vskip 0.1in \setlength{\unitlength}{1.2mm}
\begin{picture}(60,30)
\put(0, 20){\circle*{1}}
 \put(10, 20){\circle*{1}}
\put(5,25){\circle*{1}} \put(0, 30){\circle*{1}}
 \put(10,30){\circle*{1}}
 \put(0, 20){\line(1, 0){10.0}}
 \put(0, 20){\line(1, 1){5.0}}
 \put(10, 20){\line(-1, 1){5.0}}
 \put(0, 30){\line(1, 0){10.0}}
 \put(0, 30){\line(1, -1){5.0}}
 \put(10, 30){\line(-1, -1){5.0}}
 \put(0, 20){\line(0, 1){10.0}}
 \put(10, 30){\line(0, -1){10.0}}
 \put(20, 20){\circle*{1}}
 \put(30, 20){\circle*{1}}
\put(40,20){\circle*{1}}
 \put(20, 30){\circle*{1}}
 \put(30,30){\circle*{1}}
 \put(40, 30){\circle*{1}}
 \put(50,30){\circle*{1}}
\put(50, 20){\circle*{1}}
 \put(60,30){\circle*{1}}
 \put(70, 30){\circle*{1}}
 \put(70, 20){\circle*{1}}
\put(80,30){\circle*{1}}
 \put(80, 20){\circle*{1}}
\put(80, 30){\line(0, -1){10.0}} \put(70, 30){\line(0, -1){10.0}}
\multiput(73,25)(2,0){3}{\circle*{0.3}{3}}
\put(20,30){\line(1,0){10.0}} \put(30, 30){\line(1, -1){10.0}}
\put(50,30){\line(1,0){10.0}}
   \put(40, 30){\line(1, -1){10.0}}
   \put(20, 30){\line(1, -1){10.0}}
    \put(20, 20){\line(1, 0){10.0}}
\put(40, 30){\line(1, 0){10.0}} \put(30, 20){\line(1, 0){10.0}}
      \put(20, 30){\line(0, -1){10.0}}
      \put(20, 20){\line(3, 1){30.0}}

\multiput(31,30)(3,0){3}{\line(1,0){2}}

      \multiput(41,20)(3,0){3}{\line(1,0){2}}
\put(20,10){\makebox(8,1)[l]{$(4, 3^{a},2,1^{n-2-a})$}}
\put(20,5){\makebox(8, 1)[l]{Figure 25}}
\end{picture}
\vskip 0.1in \setlength{\unitlength}{1.2mm}
\begin{picture}(60,30)
\put(0, 20){\circle*{1}}
 \put(10, 20){\circle*{1}}
\put(5,25){\circle*{1}} \put(0, 30){\circle*{1}}
 \put(10,30){\circle*{1}}
 \put(0, 20){\line(1, 0){10.0}}
 \put(0, 20){\line(1, 1){5.0}}
 \put(10, 20){\line(-1, 1){5.0}}
 \put(0, 30){\line(1, 0){10.0}}
 \put(0, 30){\line(1, -1){5.0}}
 \put(10, 30){\line(-1, -1){5.0}}
 \put(0, 20){\line(0, 1){10.0}}
 \put(10, 30){\line(0, -1){10.0}}
 \put(20, 20){\circle*{1}}
 \put(30, 20){\circle*{1}}
\put(40,20){\circle*{1}}
 \put(20, 30){\circle*{1}}
 \put(30,30){\circle*{1}}
 \put(40, 30){\circle*{1}}
 \put(50,30){\circle*{1}}
\put(50, 20){\circle*{1}}
 \put(60,30){\circle*{1}}
 \put(70, 30){\circle*{1}}
\put(80,30){\circle*{1}}
 \put(90, 30){\circle*{1}}
 \put(80, 20){\circle*{1}}
 \put(90,20){\circle*{1}}
\put(80, 30){\line(0, -1){10.0}}
 \put(90, 30){\line(0, -1){10.0}}
\multiput(83,25)(2,0){3}{\circle*{0.3}{3}}
\put(20,30){\line(1,0){10.0}} \put(30, 30){\line(1, -1){10.0}}
\put(50,30){\line(1,0){10.0}}
   \put(40, 30){\line(1, -1){10.0}}
   \put(20, 30){\line(1, -1){10.0}}
    \put(20, 20){\line(1, 0){10.0}}
\put(40, 30){\line(1, 0){10.0}} \put(30, 20){\line(1, 0){10.0}}
      \put(20, 30){\line(0, -1){10.0}}
      \put(20, 20){\line(3, 1){30.0}}

\multiput(31,30)(3,0){3}{\line(1,0){2}}
      \multiput(61,30)(3,0){3}{\line(1,0){2}}
      \multiput(41,20)(3,0){3}{\line(1,0){2}}
\put(20,10){\makebox(8,1)[l]{$(4, 3^{a},2^b,1^{n-1-a-b})$}}
\put(20,5){\makebox(8, 1)[l]{Figure 26}}
\end{picture}
\vskip 0.1in \setlength{\unitlength}{1.2mm}
\begin{picture}(60,30)
\put(0, 20){\circle*{1}}
 \put(10, 20){\circle*{1}}
\put(5,25){\circle*{1}} \put(0, 30){\circle*{1}}
 \put(10,30){\circle*{1}}
 \put(0, 20){\line(1, 0){10.0}}
 \put(0, 20){\line(1, 1){5.0}}
 \put(10, 20){\line(-1, 1){5.0}}
 \put(0, 30){\line(1, 0){10.0}}
 \put(0, 30){\line(1, -1){5.0}}
 \put(10, 30){\line(-1, -1){5.0}}
 \put(0, 20){\line(0, 1){10.0}}
 \put(10, 30){\line(0, -1){10.0}}
\put(20,20){\circle*{1}}
 \put(20,30){\circle*{1}}
 \put(30,20){\circle*{1}}
 \put(40, 20){\circle*{1}}
 \put(30,30){\circle*{1}}
 \put(40, 30){\circle*{1}}
 \put(50,30){\circle*{1}}
 \put(50, 20){\circle*{1}}
 \put(60,30){\circle*{1}}
 \put(60, 20){\circle*{1}}
\put(70,30){\circle*{1}} \put(70,20){\circle*{1}} \put(70,
30){\line(0, -1){10.0}}
 \put(50, 30){\line(0, -1){10.0}}
 \put(30, 30){\line(0, -1){10.0}}
\put(60, 30){\line(0, -1){10.0}}
 \put(30, 30){\line(1, 0){10.0}}

\put(20, 30){\line(1, 0){10.0}} \put(20, 20){\line(1, 0){10.0}}
\put(30, 20){\line(1, 0){10.0}}
\multiput(63,25)(2,0){3}{\circle*{0.3}{3}}
\put(20,10){\makebox(8,1)[l]{$(4, 3^{6},1^{n-7})$}}
\put(20,5){\makebox(8, 1)[l]{Figure 27}}

\end{picture}
\vskip 0.1in \setlength{\unitlength}{1.2mm}
\begin{picture}(60,30)
\put(0, 20){\circle*{1}}
 \put(10, 20){\circle*{1}}
\put(5,25){\circle*{1}} \put(0, 30){\circle*{1}}
 \put(10,30){\circle*{1}}
 \put(0, 20){\line(1, 0){10.0}}
 \put(0, 20){\line(1, 1){5.0}}
 \put(10, 20){\line(-1, 1){5.0}}
 \put(0, 30){\line(1, 0){10.0}}
 \put(0, 30){\line(1, -1){5.0}}
 \put(10, 30){\line(-1, -1){5.0}}
 \put(0, 20){\line(0, 1){10.0}}
 \put(10, 30){\line(0, -1){10.0}}
\put(20,20){\circle*{1}}
 \put(20,30){\circle*{1}}
 \put(30,20){\circle*{1}}
 \put(40, 20){\circle*{1}}
 \put(30,30){\circle*{1}}
 \put(40, 30){\circle*{1}}
 \put(50,30){\circle*{1}}
 \put(50, 20){\circle*{1}}
 \put(60,30){\circle*{1}}
 \put(60, 20){\circle*{1}}
\put(70,30){\circle*{1}}
 \put(70, 20){\circle*{1}}
 \put(50, 30){\line(0, -1){10.0}}
 \put(30, 30){\line(0, -1){10.0}}
 \put(20, 30){\line(0, -1){10.0}}
 \put(40, 30){\line(0, -1){10.0}}
 \put(20, 30){\line(3, -1){30.0}}
 \put(20, 20){\line(3, 1){30.0}}
\put(60, 30){\line(0, -1){10.0}}
 \put(30, 30){\line(1, 0){10.0}}
 \put(70, 30){\line(0, -1){10.0}}
\put(20, 30){\line(1, 0){10.0}} \put(20, 20){\line(1, 0){10.0}}
\put(30, 20){\line(1, 0){10.0}}
\multiput(41,30)(3,0){3}{\line(1,0){2}}
\multiput(41,20)(3,0){3}{\line(1,0){2}}
\multiput(63,25)(2,0){3}{\circle*{0.3}{3}}
\put(20,10){\makebox(8,1)[l]{$(4, 3^{a},1^{n-1-a})$}}
\put(20,5){\makebox(8, 1)[l]{Figure 28}}
\end{picture}
\vskip 0.1in \setlength{\unitlength}{1.2mm}
\begin{picture}(60,30)
\put(0, 20){\circle*{1}}
 \put(10, 20){\circle*{1}}
\put(5,25){\circle*{1}} \put(0, 30){\circle*{1}}
 \put(10,30){\circle*{1}}
 \put(0, 20){\line(1, 0){10.0}}
 \put(0, 20){\line(1, 1){5.0}}
 \put(10, 20){\line(-1, 1){5.0}}
 \put(0, 30){\line(1, 0){10.0}}
 \put(0, 30){\line(1, -1){5.0}}
 \put(10, 30){\line(-1, -1){5.0}}
 \put(0, 20){\line(0, 1){10.0}}
 \put(10, 30){\line(0, -1){10.0}}
\put(20,20){\circle*{1}}
 \put(20,30){\circle*{1}}
 \put(30,20){\circle*{1}}
 \put(40, 20){\circle*{1}}
 \put(30,30){\circle*{1}}
 \put(40, 30){\circle*{1}}
 \put(50,30){\circle*{1}}
 \put(50, 20){\circle*{1}}
 \put(60,30){\circle*{1}}
 \put(60, 20){\circle*{1}}
\put(40, 30){\line(0, -1){10.0}}
 \put(50, 30){\line(0, -1){10.0}}
 \put(20, 30){\line(0, -1){10.0}}
\put(60, 30){\line(0, -1){10.0}}
 \put(20, 30){\line(1, -1){10.0}}
\put(20, 30){\line(1, 0){10.0}}
\multiput(53,25)(2,0){3}{\circle*{0.3}{3}}
\put(20,10){\makebox(8,1)[l]{$(4, 3^{5},1^{n-6})$}}
\put(20,5){\makebox(8, 1)[l]{Figure 29}}
\end{picture}
\vskip 0.1in \setlength{\unitlength}{1.2mm}
\begin{picture}(60,30)
\put(0, 20){\circle*{1}}
 \put(10, 20){\circle*{1}}
\put(5,25){\circle*{1}} \put(0, 30){\circle*{1}}
 \put(10,30){\circle*{1}}
 \put(0, 20){\line(1, 0){10.0}}
 \put(0, 20){\line(1, 1){5.0}}
 \put(10, 20){\line(-1, 1){5.0}}
 \put(0, 30){\line(1, 0){10.0}}
 \put(0, 30){\line(1, -1){5.0}}
 \put(10, 30){\line(-1, -1){5.0}}
 \put(0, 20){\line(0, 1){10.0}}
 \put(10, 30){\line(0, -1){10.0}}
\put(20,20){\circle*{1}}

 \put(30,20){\circle*{1}}
 \put(40, 20){\circle*{1}}
 \put(30,30){\circle*{1}}
 \put(40, 30){\circle*{1}}

 \put(50, 20){\circle*{1}}
 \put(60,30){\circle*{1}}
 \put(60, 20){\circle*{1}}
\put(70,30){\circle*{1}}
 \put(70, 20){\circle*{1}}

 \put(30, 30){\line(0, -1){10.0}}
  \put(30, 30){\line(1, -1){10.0}}

\put(60, 30){\line(0, -1){10.0}}
 \put(30, 30){\line(1, 0){10.0}}
 \put(40, 20){\line(1, 0){10.0}}
 \put(70, 30){\line(0, -1){10.0}}
 \put(20, 20){\line(1, 0){10.0}}
\put(30, 20){\line(1, 0){10.0}}

\multiput(63,25)(2,0){3}{\circle*{0.3}{3}}
\put(20,10){\makebox(8,1)[l]{$(4, 3^{7},1^{n-8})$}}
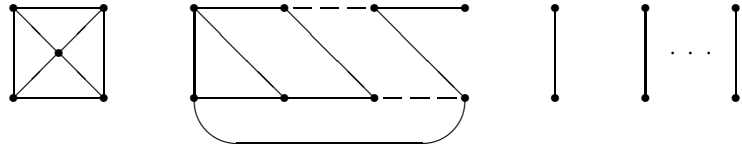
\put(20,5){\makebox(8, 1)[l]{Figure 30}}
\end{picture}
\vskip 0.1in \setlength{\unitlength}{1.2mm}
\begin{picture}(60,30)
\put(0, 20){\circle*{1}}
 \put(10, 20){\circle*{1}}
\put(5,25){\circle*{1}} \put(0, 30){\circle*{1}}
 \put(10,30){\circle*{1}}
 \put(0, 20){\line(1, 0){10.0}}
 \put(0, 20){\line(1, 1){5.0}}
 \put(10, 20){\line(-1, 1){5.0}}
 \put(0, 30){\line(1, 0){10.0}}
 \put(0, 30){\line(1, -1){5.0}}
 \put(10, 30){\line(-1, -1){5.0}}
 \put(0, 20){\line(0, 1){10.0}}
 \put(10, 30){\line(0, -1){10.0}}
 \put(20, 20){\circle*{1}}
 \put(30, 20){\circle*{1}}
\put(40,20){\circle*{1}}
 \put(20, 30){\circle*{1}}
 \put(30,30){\circle*{1}}
 \put(40, 30){\circle*{1}}
 \put(50,30){\circle*{1}}
\put(50, 20){\circle*{1}}
 \put(60,30){\circle*{1}}
 \put(70, 30){\circle*{1}}
 \put(70, 20){\circle*{1}}
 \put(80,30){\circle*{1}}
 \put(60, 20){\circle*{1}}
 \put(80, 20){\circle*{1}}

 \put(70, 30){\line(0, -1){10.0}}
 \put(80, 30){\line(0, -1){10.0}}
 \put(60, 30){\line(0, -1){10.0}}
\multiput(73,25)(2,0){3}{\circle*{0.3}{3}}
\put(20,30){\line(1,0){10.0}} \put(30, 30){\line(1, -1){10.0}}

   \put(40, 30){\line(1, -1){10.0}}
   \put(20, 30){\line(1, -1){10.0}}
    \put(20, 20){\line(1, 0){10.0}}
\put(40, 30){\line(1, 0){10.0}} \put(30, 20){\line(1, 0){10.0}}
      \put(20, 30){\line(0, -1){10.0}}

\multiput(31,30)(3,0){3}{\line(1,0){2}}
 \put(35, 20){\oval(30,10)[b]}

      \multiput(41,20)(3,0){3}{\line(1,0){2}}
\put(20,10){\makebox(8,1)[l]{$(4, 3^{a},1^{n-1-a})$}}
\put(20,5){\makebox(8, 1)[l]{Figure 31}}

\end{picture}
\end{document}